\documentstyle[11pt, sec]{article}
\newfont{\Bbb}{msbm10 scaled \magstep1}

\newcommand\uc{\underline{\hbox{\Bbb C}}}

\newcommand\bC{\hbox{\Bbb C}}

\newcommand\bR{\hbox{\Bbb R}}
\newcommand\bZ{\hbox{\Bbb Z}}

\newcommand\bS{\hbox{\Bbb S}}

\newfont{\es}{eusm10 scaled \magstep1}
\newfont{\ses}{eufm8 scaled \magstep1}
\newfont{\gt}{eufb10 scaled \magstep1}
\newfont{\sg}{eufb8 scaled \magstep1}
\newfont{\goth}{eufb10 scaled \magstep2}

\newcommand{\re}{\hbox{\gt Re}}

\newcommand{\gc}{\hbox{\gt c}}

\newcommand{\gf}{\hbox{\gt f}}
\newcommand{\gG}{\hbox{\gt G}}

\newcommand{\dir}{\hbox{\es D}}
\newcommand{\en}{\hbox{\es E}}

\newcommand{\dol}{\hbox{\goth d}}
\newcommand{\can}{\hbox{\es K}}

\newcommand{\da}{\dot{a}}
\newcommand{\dps}{\dot{\psi}}

\newcommand{\de}{{\delta}}

\newcommand{\ve}{{\varepsilon}}

\newcommand{\si}{{\sigma}}

\def\L{{\cal L}}
\newcommand{\s}{{\cal S}}

\def\ra{\rightarrow}
\def\hra{\hookrightarrow}
\def\be{\begin{equation}}
\def\ee{\end{equation}}
\def\lan{\langle}
\def\ran{\rangle}

\newcommand{\naf}{{}^\flat\nabla}
\newcommand{\onaf}{{}^\flat\overline{\nabla}}

\newtheorem{theorem}{Theorem}[section]
\newtheorem{lemma}[theorem]{Lemma}

\newtheorem{remark}[theorem]{Remark}

\begin{document}

\title{Seiberg-Witten Equations on  Tubes}
\author{Liviu I. Nicolaescu\thanks{{\bf Current address}: Dept.of Math.,University of Michigan, Ann Arbor, MI 48109-1109, USA;  liviu@math.lsa.umich.edu}}
\date{December, 1996}
\maketitle

\addcontentsline{toc}{section}{Introduction}

\begin{center}
{\bf Introduction}
\end{center}

\bigskip

  In \cite{Ni} we began   a study of the 3-dimensional Seiberg-Witten  equations
   on Seifert manifolds with  two goals in mind: to compute the Seiberg-Witten-Floer (SWF) homology of these manifolds and  ultimately, to produce gluing formul{\ae} for the 4-dimensional Seiberg-Witten invariants.

A first difficulty one must overcome    has to do with the less than obvious nature 
 of the solutions of the 3-dimensional Seiberg-Witten equations. We  dealt with 
 this issue in \cite{Ni} where  we studied  the behavior of these solutions as  
 the Seifert fibration collapses onto its base (i.e. the  background metric is shrunk
  in the fiber direction). As the metric is deformed, the solutions of the SW equations converge to solutions of  some {\em adiabatic Seiberg-Witten equations}. These  are variational equations and can be solved {\em explicitly}.  Moreover, these adiabatic equations  are very simple {\em zeroth order perturbations} of the original ones  which suggests that the Morse theory for the adiabatic equations produces the same results as the original ones (which may have to be perturbed anyway to be placed in a generic framework).  

A key fact established in \cite{Ni} was that,  in the case of a smooth $S^1$-bundle $N$  over a Riemann surface $\Sigma$  equipped with a product-like metric with sufficiently short fibers,  the adiabatic Morse function is Bott nondegenerate along the irreducible part of its critical set.  This fact makes the adiabatic theory even more tempting to use for Floer theory  computations.  

The Bott  extension of Morse theory (in the form described  for example in
\cite{AB}) describes a spectral sequence associated to a Morse-Bott function converging to
 the cohomology of the background manifold. This approach can be extended to
 our infinite dimensional situation as in \cite{Fu} for the instanton homology.

The (co)boundary operators of this spectral sequence are defined in terms of  the 
tunnelings between different components of the critical set,
i.e. connecting trajectories of the gradient flow. In the case at hand, the  gradient
 flow equations are in fact  (zeroth order) perturbations of the 4-dimensional Seiberg-Witten equations on the tube ${\bR}\times N$.   

The main results of this paper (Theorem \ref{th: vanish}) will show that these 
 tunnelings {\em do not exist} provided the fibers of $N$ are sufficiently short.  This  implies the boundary operators are trivial and 
  thus, the adiabatic theory leads to a {\em perfect} Morse function. To actually compute the (co)homology (as a {\em graded} object) one has to  compute several spectral flows. We will address this issue elsewhere.

The proof of Theorem \ref{th: vanish} is conceptually very simple.    First of all, the results of \cite{Ni} show that the critical points of  the  adiabatic Morse functional    do not change as   we shrink the fibers of $N$. In particular, there exists a {\em positive} {\em lower bound} for  the $L^\infty$ norms of these solutions and this bound is {\em independent of the shrinking geometry}.  We next show that  (when they exist) 
the    tunnelings  can be used to  produce an   {\em effective upper bound} for  the $L^\infty$ norms of the above critical points. This upper bound {\em converges to zero} as the fibers of $N$ become shorter and shorter. Thus tunnelings cannot exist if $N$ has short fibers.  The key ingredient in the proof of this effective $L^\infty$-estimate is  M\"{o}ser's iteration technique   in which we carefully keep track of the dependence of the best Sobolev embedding constants on  the  shrinking geometry.

\noindent {\bf Note} After this paper was completed we learned  of the paper \cite{MOY}
where   these tunnelings are studied   via algebraic geometric techniques. 
It is however not clear  whether those  results imply the adiabatic  dissapearance of
tunnelings  proved in this paper.

\tableofcontents

\section{Seiberg-Witten equations on smooth circle bundles}

In this section we briefly  survey some basic facts established in \cite{Ni}
(see also \cite{MOY}) about the 3-dimensional Seiberg-Witten equations on smooth circle fibrations and then we introduce the 4-dimensional Seiberg-Witten equations on tubes.   

\subsection{The differential geometric background}
 
Consider $\ell \in {\bZ}$ and  denote by $N=N_\ell$ the total space of a degree
$\ell$ principal $S^1$ bundle over a  compact oriented surface of genus $g$: $S^1 \hookrightarrow N_\ell \stackrel{\pi}{\ra} \Sigma$.  Denote by $\zeta\in{\rm Vect}\,(N)$ the infinitesimal generator of the $S^1$ action.  $N$ has a natural orientation  which can be described using any   any splitting $TN=\lan \zeta\ran \oplus \pi^*T\Sigma$ determined by an arbitrary connection.

Assume $\Sigma$ is equipped with a Riemann metric $h_b$ such that ${\rm vol}_{h_b}(\Sigma)= \pi$. Pick a connection form ${\bf i}\eta\in {\bf i}\Omega^1(N)$ such that
 \[
 -d\eta= 2 \ell dv_{h_b}.
 \]
This choice is possible since $\frac{-1}{2\pi}d\eta$ represents the first Chern class of $N$ which is $\ell$.  For each $\delta \geq 1$ define a metric $h_\delta$ on $N$ by
\[
h_\delta=\delta^{-2}\eta\otimes \eta \oplus \pi^* h_b.
\]
When $\delta =1$ we will write $h$ instead of $h_1$.

Using this metric we can orthogonally  split $T^*N\cong \lan \eta \ran \oplus \pi^*T^*\Sigma$  and  this defines in a natural way an orientation on $N$.  If $\ast_\delta$ denotes the Hodge $\ast$ operator of the metric $h_\delta$ we get 
\[
d\eta_\delta =2\lambda_\delta \ast_\delta \eta_\delta
\]
where $\eta_\delta=\delta^{-1}\eta$ and $\lambda_\delta=-\ell\delta^{-1}$. Again we set $\lambda=\lambda_1$.

Fix  a local, oriented $h$-orthonormal  coframe $\eta_0=\eta, \eta_1, \eta_2$ on $N$ and denote by $\zeta^0=\zeta,\zeta^1,\zeta^2$ its dual frame. The bundle $\lan \eta\ran^\perp$ has a natural complex structure  locally defined by the correspondences $\eta_1\mapsto-\eta_2 \mapsto-\eta_1$.  In this way we get a complex line bundle ${\can}\ra N$. It is isomorphic with the pullback of the canonical line bundle $K_\Sigma$ of the base.    We have a   splitting
\be
T^*N\otimes {\bC}\cong {\uc}\oplus {\can}\oplus {\can}^{-1}
\label{eq: can}
\ee
where ${\uc}$ denotes  the trivial complex line bundle (over $N$). Note that the complex hermitian line bundle ${\can}$ has {\em two} natural connections on it. A connection $\nabla=\nabla^h$ induced by the  Levi-Civita connection of the metric $h$ via the decomposition (\ref{eq: can}) and a connection $\nabla^\infty$ induced by pullback from the Levi-Civita connection on $K_\Sigma$.

For any complex vector bundle $E\ra N$ we get
\be
T^*N\otimes E\cong E\, \oplus\,  {\can}\otimes E\, \oplus \, {\can}^{-1}\otimes E.
\label{eq: conn}
\ee
Thus,  any connection $\nabla:C^\infty(E)\ra C^\infty(T^*N\otimes E)$ naturally splits in three parts

\noindent $\bullet$ $\nabla_\zeta:C^\infty(E)\ra C^\infty(E)$.

\noindent $\bullet$ ${}^\flat\nabla={\ve}\otimes (\nabla_1-{\bf i}\nabla_2):C^\infty(E)\ra C^\infty({\can}\otimes E)$ 

\noindent$\bullet$  ${}^\flat\bar{\nabla}=\bar{\ve}\otimes (\nabla_1+{\bf i}\nabla_2):C^\infty(E)\ra C^\infty({\can}^{-1}\otimes E)$.

\noindent where ${\ve}=2^{-1/2}(\eta_1+{\bf i}\eta_2)$, $\bar{\ve}=2^{-1/2}(\eta_1-{\bf i}\eta_2)$ and $\nabla_j=\nabla_{\zeta^j}$ for $j=0,1,2$. (Note that 
$\nabla=\eta\otimes \nabla_0 +2^{-1/2}({}^\flat\nabla +{}^\flat\bar{\nabla})$.

Since $TN$ is trivial, $N$ admits $spin$ structures and  they are parameterized by the square roots of  ${\can}$. We pick a square root of ${\can}$ as a pullback of a fixed square root $K^{1/2}$ of $K_\Sigma$.  Once this $spin$ structure is fixed, the $spin^c$ structures are parameterized bijectively by the family of isomorphism classes of complex line bundles on  $N$. 

The bundle of complex  spinors ${\bS}_\sigma$ corresponding to the $spin^c$  structure $\sigma$  determined by a line bundle  ${\cal L}={\cal L}_\sigma \ra N$   is
\be
{\bS}={\bS}_\sigma={\can}^{-1/2}\otimes {\cal L}\; \oplus \; {\can}^{1/2}\otimes {\cal L}.
\label{eq: spinor}
\ee
If $A$ is a connection on ${\cal L}$  then we can tensor it either with $\nabla^h$, or with $\nabla^\infty$ to obtain two pairs of connections ($\nabla^A, \nabla^{A,\infty})$,  a pair on    ${\can}^{-1/2}\otimes {\cal L}$ and  the other on ${\can}^{1/2}\otimes {\cal L}$.  We can now construct the   operators
\[ 
\nabla^A_\zeta, \nabla_\zeta^{A,\infty}: C^{\infty}({\can}^{\pm 1/2}\otimes {\cal L})  \ra  C^{\infty}({\can}^{\pm 1/2}\otimes {\cal L} ) 
\]
\[
{\naf}^A={\naf}^{A,\infty}: C^\infty({\can}^{-1/2}\otimes {\cal L}) \ra C^\infty({\can}^{1/2}\otimes {\cal L})
\]
and 
\[
{\onaf}^A={\onaf}^{A,\infty}:C^\infty({\can}^{1/2}\otimes {\cal L}) \ra C^\infty({\can}^{-1/2}\otimes {\cal L}) 
\]
Now form the operators $Z_A, Z_{A,\infty}, T_A:C^\infty({\bS})\ra C^\infty({\bS})$ which in terms of the splitting  (\ref{eq: spinor}) have the block decompositions
\[
Z_A= \left[
\begin{array}{cc}
{\bf i}\nabla^A_\zeta & 0 \\
0 & -{\bf i}\nabla^A_\zeta 
\end{array}
\right] \;\; Z_{A,\infty}=\left[
\begin{array}{cc}
{\bf i}\nabla^{A,\infty}_\zeta & 0 \\
0 & -{\bf i}\nabla^{A,\infty}_\zeta 
\end{array}
\right] 
\]
\[
T_A=\left[ 
\begin{array}{cc}
0 & {\onaf}^A \\
{\naf}^A & 0 
\end{array}
\right]
\]
These operators are formally selfadjoint and satisfy the following fundamental  identities.
\be
Z_{A,\infty}=Z_A+\lambda/2 =Z_A-\ell/2.
\label{eq: za}
\ee
\be
\{Z_A, T_A\}:=Z_AT_A+T_AZ_A=-\lambda T_A + {\bf i}\left[
\begin{array}{cc}
0 & F_A^{0,1} \\
F_A^{1,0} & 0 
\end{array}
\right]
\label{eq: zt}
\ee
where $F_A$ denotes the curvature of $A$, $F_A^{0,1}=\bar{\ve}\otimes F_A(\zeta, \zeta_1+{\bf i} \zeta_2)$ and $F_A^{1,0}={\ve}\otimes F_A(\zeta,\zeta_1-{\bf i}\zeta_2)$.
The last two equalities imply
\be
\{Z_{A,\infty},T_A\}={\bf i}\left[
\begin{array}{cc}
0 & F_A^{0,1} \\
F_A^{1,0} & 0 
\end{array}
\right]
\label{eq: azt}
\ee
The Dirac operator on ${\bS}$ defined by the connection $A$ is 
\be
{\dir}_A= Z_A + T_A +\lambda=Z_A+ T_A -\ell
\label{eq: dir}
\ee
Define the {\em adiabatic  Dirac operator} by
\be
{\bf D}_A =Z_{A,\infty}+T_A= {\dir}_A-\lambda/2.
\label{eq: adir}
\ee
The above  constructions depend on  the parameter $\delta$. For the metric $h_\delta$  one defines $Z_{A,\infty, \delta}=\delta Z_{A,\infty}$, $Z_{A,\delta}=Z_{A,\infty, \delta}-\lambda_\delta/2$, ${\dir}_{A,\delta}=Z_{A,\delta}+T_A+\lambda_\delta$ and ${\bf D}_{A,\delta}=Z_{A,\infty, \delta}+T_A$.

\subsection{The 3-dimensional Seiberg-Witten equations}
The data entering the Seiberg-Witten equations are the  following.

\noindent {\bf  (a)} A  $spin^c$ structure $\sigma$ on $N$ determined by the line 
bundle $\L={\L}_\sigma$.

\noindent {\bf (b)} A connection $A$ of ${\L}\ra N$.

\noindent {\bf (c)} A    spinor $\psi$ i.e. a section of the complex spinor    
bundle ${\bS}_{\sigma}$ associated to the given  ${\rm spin}^c$ structure $\sigma$.

The connection $A$  defines a geometric Dirac operator  ${\dir}_A$  on
${\bS}_{\sigma}$. The Seiberg-Witten equations are
\[
(SW):\;\;\left\{
\begin{array}{rcr}
{\dir}_{A}\psi  & = & 0 \\
{\bf c}(\ast F_{A})& =  & \tau (\phi)
\end{array}
\right.
\]
where $\ast$ is the Hodge $\ast$-operator of the metric $g$, $\tau(\psi)=\bar{\psi}\otimes \psi-\frac{1}{2}|\psi|^2\in {\rm End}\,({\bS})$  and  
${\bf c}$ is the Clifford multiplication  defined in terms of the background metric $h$.  If the metric $h$ is replaced by the metric $h_\delta$ we get new equations $(SW)_\delta$  obtainable from $(SW)$ via the substitutions
\[
{\dir}_A\mapsto {\dir}_{A,\delta},\;\;\ast\mapsto \ast_\delta,\;\;\;{\bf c}\mapsto {\bf c}_\delta.
\] 
The {\em adiabatic  Seiberg-Witten equations} are
\[
(SW)_\infty:\;\;\left\{
\begin{array}{rcr}
{\bf D}_{A}\psi  & = & 0 \\
{\bf c}(\ast F_{A})& =  & \tau (\phi)
\end{array}
\right.
\]
For each $\delta\geq 1$ we get an adiabatic equation $(SW)_{\infty, \delta}$  derived from  $(SW)_\infty$  in the same way $(SW)_\delta$ is obtained from $(SW)$.

Using the decomposition ${\bS}={\can}^{-1/2}\otimes {\L}\; \oplus \; {\can}^{1/2}\otimes {\L}$ we can decompose any spinor $\psi$ as $\psi =\alpha \oplus \beta$  and  in terms of this decomposition the Seiberg-Witten  equations  can be rephrased as
\be
\left\{
\begin{array}{rrrcl}
{\bf i}{\nabla}_\zeta^A \alpha & +  {\onaf}^A \beta & +\lambda  \alpha &  = & 0 \\
 & & & & \\
{\naf}^A \alpha & - {\bf i}{\nabla}^A_\zeta \beta  & +\lambda  \beta & = &0 \\
&&&&\\

& &  \frac{1}{2}(|\alpha|^2-|\beta|^2) & = &{\bf i}F_A(\zeta_1, \zeta_2) \\
& & & & \\
& & {\bf i}\alpha \bar{\beta}&= &\bar{\ve}\otimes F_A(\zeta_1
+{\bf i} \zeta_2, \zeta) =-F_A^{0,1}
\end{array}
\right.
\label{eq: sw}
\ee
The adiabatic Seiberg-Witten equations  are obtained from the above  by replacing $\lambda$ with $\lambda/2$. 

The Seiberg-Witten equations have a variational nature. Fix a smooth connection
$A_0$ on ${\cal L}$ and define
\[
{\gf}: L^{1,2}({\bS}\oplus {\bf i}T^*N) \ra {\bR}
\]
by
\[
{\gf}(\psi, a)=\frac{1}{2}\int_Na\wedge(F_{A_0}+ F_{A_0+a}) +\frac{1}{2}\int_N  \lan \psi,
{\dir}_{A_0+a}\psi \ran dv_h.
\]
The differential of ${\gf}$ at a point ${\gc}=(\phi,a)$ is
\[
d_{{\gc}}{\gf}(\dot{\phi}, \dot{a})=\int_N\lan \dot{a}, {\bf c}^{-1}(\tau(\phi))-\ast F_{A_0+a}\ran
dv_g +\int_N{\re}\lan \dot{\phi}, {\dir}_{A_0+a}\phi\ran dv_g.
\]
Hence, the solutions of $(SW)$ are critical points of ${\gf}$. Similarly, the solutions of $(SW)_\infty$ are  critical points of 
\[
{\en}(\psi,a)= \frac{1}{2}\int_Na\wedge(F_{A_0}+ F_{A_0+a}) +\frac{1}{2}\int_N  \lan \psi,
{\bf D}_{A_0+a}\psi \ran dv_h
\]
The solutions of $(SW)_{\infty, \delta}$ are critical points of a functional
${\en}_\delta$, defined as ${\en}$ but in terms of the metric $h_\delta$.

Obviously, the adiabatic equations are zeroth order perturbations of the original ones.  It is less obvious  though,  that the solutions of $(SW)_\delta$  converge as $\delta\ra \infty$ to the solutions of $(SW)_\infty$ (see \cite{Ni}).   The remarkable feature of the adiabatic equations  is that {\em they can be solved explicitly}, due mainly to the identity (\ref{eq: azt}).   In fact, for each $\delta$  we  get a space of solutions of $(SW)_{\infty,\delta}$ but remarkably, these spaces turn out to be {\em independent of} $\delta$!

To present the explicit description of the solutions of $(SW)_\infty$ let us first identify the space of $spin^c$ structures on $N$  with $H^2(N, {\bZ})$. Using the Thom-Gysin exact sequence for the $S^1$-bundle $N\stackrel{\pi}{\ra}\Sigma$ we deduce
\[
H^2(N)\cong H^1(\Sigma, {\bZ})\oplus \pi^*H^2(\Sigma,{\bZ})\cong {\bZ}^{2g}\oplus {\bZ}_{|\ell|}.
\]
Denote by ${\cal S}_\infty={\cal S}_\infty(\L)$ the space  of solutions of $SW_\infty$ modulo the action of the gauge group ${\gG}={\rm Aut}\, (L)$. In \cite{Ni} we proved the following facts. We assume $\ell \neq 0$.

\noindent {\bf A.} If $c_1({\L})\in H^2(N, {\bZ})$ {\em is not a torsion class} then ${\cal S}_\infty =\emptyset$.

\medskip

\noindent {\bf B.} If $c_1({\L})$ is a torsion class, $c_1=\kappa \in {\bZ}_{|\ell|}$,  then ${\en}$ is ${\gG}$-{\em invariant} (!)  and  the  adiabatic  moduli space ${\cal S}_\infty$  consists of the following.

\noindent $\bullet$ A reducible part ${\cal S}_\infty^{red}$  which is homeomorphic with a torus $T^{2g}$.   A pair  $(\psi, A)$ lies in  ${\cal S}_\infty^{red}$ if and only if $\psi =0$ and $A$ is a flat  connection on ${\L}$ with holonomy along the fibers $\exp(2\pi {\bf i}\kappa/\ell)$.

\noindent $\bullet$ An irreducible part ${\cal S}^{irr}$.   Set  $I_\kappa=\{  n\in {\bZ}\; ; \; 0<|n|\leq g-1,\;\;n\equiv \kappa\; {\rm mod}\; \ell\}$.   
 $(\alpha \oplus \beta, A)$ lies in ${\cal S}_\infty^{irr}$ if and only if the following happen. 

(i)  There exists a complex line bundle $L\ra \Sigma$  and a connection $A_\Sigma$ on $L$ such that $\deg L\in I_\kappa$, $\pi^*L={\L}$ and 
$\pi^*A_\Sigma=A$.

(ii)  If $\deg L < 0$ then  $\alpha =0$ and $\beta$ is the pullback of a holomorphic  section of $K^{1/2}\otimes L$. Moreover,
\be
\frac{1}{4\pi}\|\beta\|^2 =-\deg L.
\label{eq: n1}
\ee

(iii) If $\deg L >0$ then $\beta =0$ and $\alpha$ is the pullback of an antiholomorphic section of $K^{-1/2}\otimes L$. Moreover,
\be
\frac{1}{4\pi}\|\alpha\|^2 =\deg L.
\label{eq: n2}
\ee

The above facts show that if $(\psi, A)\in {\cal S}_\infty^{irr}$ then $A$ induces on ${\cal L}$ a structure of $S^1$-equivariant line bundle and $A$ is an $S^1$-equivariant connection. As such,  it has  an  $S^1$-equivariant first Chern class $\hat{c}_1({\L},A)\in H^2_{S^1}(N, {\bR})$ defined using the equivariant Chern-Weil theory described e.g. in \cite{BGV}.  A classical theorem of H. Cartan (see \cite{Ni1}) identifies  $H^2_{S^1}(N)\cong H^2(\Sigma)\cong {\bR}$.  If  $({\L}, A)=\pi^*(L_\Sigma, A_\Sigma)$ (as in (i) above) then a simple computation shows  that via the above identification we have
\[
\hat{c}_1({\L}, A)=\deg L_\Sigma
\]
Hence, the number $\deg L_\Sigma$ is an {\em invariant} of the gauge
equivalence class of the solution $(\psi, A)$ of the adiabatic SW equations. We will  denote this number by $\hat{d}(\psi, A)$. We get a continuous map
\be
\hat{d}: {\cal S}_\infty^{irr}\ra I_\kappa.
\label{eq: inv}
\ee
The fibers of this map are {\em connected}  spaces, each    homeomorphic to a symmetric product  of a certain number of copies of $\Sigma$. More precisely, for $g\geq 2$ define
\[
\nu :I_\kappa \ra {\bZ}_+
\]
by
\[ 
\nu (t)= (g-1)-|t|.
\]
The space $\hat{d}^{-1}(t)$ is homeomorphic with a symmetric product of $\nu(t)$ copies of $\Sigma$.   These spaces are  spaces of {\em abelian vortices} on $\Sigma$.  Note that ${\s}_\infty^{irr}=\emptyset$ for $g=0,1$.

As we have mentioned, for each $\delta$ we get  different adiabatic equations $(SW)_{\infty, \delta}$ leading to the same moduli spaces. The energy functional ${\en}_\delta$ changes with $\delta$ and in \cite{Ni} we showed the following.

\begin{theorem}{\rm  For $\delta$ sufficiently large, ${\cal S}_\infty^{irr}$ is {\em  the nondegenerate  critical set} of ${\en}_\delta$ (modulo ${\gG}$).}
\label{th: bott}
\end{theorem}

The goal of this paper is to understand the $L^2$-gradient flow of ${\en}_\delta$ for $\delta \gg1$.

\subsection{The 4-dimensional Seiberg-Witten equation on tubes}

Consider the tube $X={\bR}\times N$ and equip it with the product orientation.  Next, fix a $spin^c$ structure $\sigma$ on
$X$ and  denote  by $\hat{\bS}_\sigma=\hat{\bS}_\sigma^+\oplus \hat{\bS}_\sigma^-$ the associated bundle of spinors.  Any connection $A$ on $\det \hat{\bS}_\sigma^+$ determines a Dirac operator 
\[
{\dol}_A: C^\infty(\hat{\bS}_\sigma^+)\ra C^\infty(\hat{\bS}_\sigma^-).
\]
The unknowns of the 4-dimensional Seiberg-Witten equations are  an even spinor  $\psi \in \Gamma(\hat{\bS}_\sigma^+)$ and a connection $A$ on $\det \hat{\bS}_\sigma^+$ such that
\[
\left\{
\begin{array}{rcl}
{\dol}_A \psi  & = & 0 \\
\frac{1}{2}\hat{\bf c}(F_A^+) & = & \tau(\psi)
\end{array}
\right.
\]
where $F_A^+= 2^{-1}(F_A+\hat{\ast}F_A)$, $\hat{\ast}$ is the Hodge $\ast$ operator on
$X$ defined by the metric $dt^2+h$, $\hat{\bf c}:\Omega^*(X)\otimes {\bC}\ra {\rm End}\,(\hat{\bS}_\sigma)$ denotes the Clifford multiplication map and 
\[
\tau(\psi)=\bar{\psi}\otimes \psi -\frac{1}{2}|\psi|^2\in {\rm End}\,(\hat{\bS}_\sigma^+)\subset  {\rm End}\, (\hat{\bS}_\sigma).
\]
The factor $1/2$  (which can be eventually absorbed by rescaling $\psi$ ) is not traditionally present in the equations but will help match the conventions of this section with those in  \S 1.2 where we  work with connections on $\det^{1/2}{\bS}_\sigma$.

Using the  special  product  structure of the tube $X$ we can further transform these equations. First of all, note that there is a bijection between the $spin^c$ structures on
$X$ and those on $N$.  Fix one such structure on $N$ and denote by ${\bS}_\sigma$ the  associated spinor bundle and by ${\bf c}$ the corresponding Clifford multiplication map. Along the slices $\{t\}\times N$ we have
\[
\hat{\bS}^\pm_\sigma\cong {\bS}_\sigma
\]
and
\be
{\bf c}(\omega)= \hat{\bf c}(dt)\hat{\bf c}(\omega)\;\;\forall \omega \in \Omega^1(N).
\label{eq: clifford}
\ee
We have seen that  the line bundle $L_\sigma=\det {\bS}_\sigma\ra N$ has a natural square root $\L_\si$.  Fix  a connection  $A_0$ on the line bundle $\L_\si$. 
Modulo  a gauge transformation, any connection on ${\bR}\times \L_\sigma$ can be
 written as $A=A_0+ a(t)$, where for each $t$, $a(t)$ is a purely imaginary 1-form 
 on $N$.  Thus $A$ can be regarded as a path of connections $A(t)=A_0+a(t)$ on $\L_\si$.
   Denote by $\hat{F}_A$ the curvature of $A$ on $X$ and, for each $t$, denote by $F_a$ the curvature of $A_0+a(t)$ on $\{t\}\times N$. 
Then
\[
\hat{F}_A= F_{a(t)}+dt\wedge \dot{a} 
\]
A simple computation shows that
\[
\hat{F}_A^+=\frac{1}{2}\left\{(F_a +\ast\dot{a}(t)) +dt\wedge (\dot{a}(t) +\ast F_a)\right\}
\]
Using the equality (\ref{eq: clifford})  we deduce that
\[
\hat{\bf c}(\hat{F}_A^+)={\bf c}(F_a+\ast \dot{a}(t))
\]
\[
{\dol}_A=\hat{\bf c}(dt)\left\{\frac{\partial}{\partial t} - {\dir}_A\right\}.
\]
A spinor $\psi\in \Gamma(\hat{\bS}_\sigma^+)$ can be regarded as a path of spinors $\psi(t)\in \Gamma({\bS}_\sigma\!\mid_{\{t\}\times N}$. Putting together all of the above  we deduce that the Seiberg-Witten equation on a tube  can be rewritten as equations for a path $A=A_0+a(t)$ of connections on $\L_\sigma$ and  a path of spinors $\psi(t)$.  More exactly, they are
\be
\left\{
\begin{array}{rcl}
\dot{\psi} & = & {\dir}_A\psi \\
{\bf c}(\dot{a}) &=& \tau(\psi) -{\bf c}(\ast F_a)
\end{array}
\right.
\label{eq: flow}
\ee
These equations describe the $L^2$ positive gradient flow of ${\gf}$.  We
can now define  {\em the adiabatic gradient flow} in an obvious fashion.
\be
\left\{
\begin{array}{rcl}
\dot{\psi} & = & {\bf D}_A\psi \\
{\bf c}(\dot{a}) &=& \tau(\psi) -{\bf c}(\ast F_a)
\end{array}
\right.
\label{eq: aflow}
\ee
We will be interested  exclusively in {\em finite energy} solutions  (or {\em
tunnelings}) of (\ref{eq:
aflow}). These are  solutions $(\psi(t), A(t))$  satisfying $\psi(t)\in
L^{1,2}_{loc}\cap L^{\infty}_{loc}$, $A(t)\in L^{1,2}_{loc}$ and
\[
\int_{-\infty}^\infty \|\dot{\psi}(t)\|^2 +\|\dot{a}(t)\|^2 dt < \infty
\]
Intuitively, the tunnelings have limits as $t\ra \pm \infty$ which
should be solutions of the 3-dimensional (adiabatic) Seiberg-Witten equations. 

The Weitzenb\"{o}ck formula for  ${\dol}_A$ has a special form on the tube.  It
splits in two parts. The first part is the Weitzenb\'{o}ck formula for ${\dir}_A$
\be
{\dir}_A^2=\nabla_A^*\nabla_A+ {\bf c}(F_A) +\frac{s}{4}
\label{eq: w1}
\ee
(where $s$ denotes the scalar curvature of  $N$) and the second part is 
\be
[\nabla_t, {\dir}_{A(t)}]:=\nabla_t {\dir}_{A(t)}-{\dir}_{A(t)}\nabla_t={\bf c}(\da)
\label{eq: w2}
\ee
where $\nabla_t$ denotes the $t$-derivative with respect to the connection
$A_0$ on ${\bR}\times \L_\si$.

\section{Tunnelings}
\setcounter{equation}{0}
In this section we study the  tunnelings of (\ref{eq: aflow}). More precisely, we will  prove the following.

\begin{theorem}{\rm  If $\delta \gg 1$ then the equation $(\ref{eq:
aflow})_\delta$ has no tunnelings.}
\label{th: vanish}
\end{theorem}

The proof will be carried out in three steps.

\noindent {\bf Step 1}\hspace{.3cm}  Produce $L^4_{loc}$ estimates  {\em
independent of} $\delta$ for the spinor component of (\ref{eq: aflow}). This 
will follow from the Weitzenb\"{o}ck formula  coupled with some {\em universal} energy
estimates.
 
\medskip

\noindent {\bf Step 2} \hspace{.3cm} Produce {\em effective} $L^\infty_{loc}$ estimates of
the spinor component in terms of the above $L^4_{loc}$ estimates and the best 
Sobolev constant of the embedding $L^{1,2}\hra L^4$ on a cylinder $[T,
T+4]\times N$. This is achieved via the M\"{o}ser's iteration technique; see
\cite{Ber} or \cite{GT}.

\medskip

\noindent {\bf Step 3}\hspace{.3cm} Conclude the proof using  the estimates in \cite{Ber}  of
 the best Sobolev constants in terms of the background geometry.

We will now supply the details.

\subsection{$L^4_{loc}$ estimates}

We will use the metric $h_\de$ on $N$  and in the  equations (\ref{eq:
aflow}) all the intervening quantities should expressed in terms of this metric.  In the
sequel, we will carefully keep track of this dependence which will be signaled by $\delta$
subscripts.

For each $T\in {\bR}$ and $L>0$ and for any tunneling $(\psi(t), a(t))$ of
(\ref{eq: aflow}) set
\[
E_\de(T, L)=E_\delta(\psi, A, T, L)= \frac{1}{2}\int_{T-L}^{T+L}dt\int_N |{\da}(t)|_\delta^2+|{\dps}(t)|^2
dv_\delta.
\]
A simple integration by parts, using  (\ref{eq: aflow}) yields the following
identity.
\be
E_\de(T,L) ={\en}_\delta(\psi(T+L), A_0 +a(T+L)) -{\en}_\delta(\psi(T-L), a(T-L)).
\label{eq: en1}
\ee

Recall the following result of \cite{MST}.

\begin{lemma}{\rm If $(\psi(t), a(t))$ is a finite energy solution of (\ref{eq:
aflow})  then there exist subsequences $t_n^\pm \ra \pm \infty$  and  solutions
$(\psi_\pm, A_\pm)\in {\cal S}_\infty(\L_\si)$ such that (modulo ${\rm Aut}\, (\L_\si)$ 
we have
\[
\lim_{t_n^\pm\ra \pm \infty} (\,\psi(t_n^\pm), A(t_n^\pm)\,)= (\psi_\pm, A_\pm)
\]
in the strong $L^{1,2}$-topology.}
\label{lemma: asymptotics}
\end{lemma}
The above lemma implies that if $c_1(\L_\si)$ is not a torsion class then  there
exist no finite energy solutions of (\ref{eq: aflow}). In the sequel, we will
assume the first Chern class is torsion. In this case the functional ${\en}$ is ${\gG}$-invariant  the following quantity is well defined and finite.
\[
G_\delta= \max\{ {\en}_\delta(\psi_+, A_+)-{\en}_\delta(\psi_-, A_-)\; ;\;
(\psi_\pm, A_\pm)\in {\s}_\infty(\L_\si)\}.
\]
$G_\delta$ is the largest energy gap between the various
components of the critical set of ${\en}_\delta$.  Using again the gauge invariance of ${\en}$,  Lemma \ref{lemma: asymptotics}  and (\ref{eq: en1})  we deduce
\be
\forall \,{\rm tunneling},\;\forall T\in {\bR},\; \forall L>0\;\;E_\delta(\psi,
A, T, L)\leq G_\delta .
\label{eq: en2}
\ee
Recall that the definition of ${\en}_\delta$ depended on a reference connection
$A_0$. We now  choose $A_0$ to be a {\em flat} connection on $\L_\si$.  Any
other connection has the form $A=A_0 +a$, $a\in \Omega^1(N)$. We can now
rewrite 
\be
{\en}_\delta=\frac{1}{2}\left\{\int_N a\wedge da + \int_N  \lan \psi,
{\bf D}_{A_0+a, \delta}\psi \ran dv_ \delta  \right\}.
\label{eq: en3}
\ee
We have the following  remarkable result.

\begin{lemma}{\rm There exists an universal constant $G$ such that}
\[
G_\delta \leq G \;\;\forall \delta \geq 1.
\]
\label{lemma: gap}
\end{lemma}

\noindent {\bf Proof}\hspace{.3cm}    The estimate follows from a  direct 
computation of the energy on the  the components of ${\cal S}_\infty$. Note
that on ${\cal S}_\infty^{red}$ we have ${\en}_\de \equiv 0$.

Along the components of ${\s}_\infty^{irr}$, the second term in the definition
of ${\en}_\delta$  vanishes (since ${\bf D}_A\psi =0$) while  the first term is 
independent of $\delta$.  More explicitly,  the energy along the component of
${\s}_\infty^{irr}$ given by the fiber $\hat{d}^{-1}(n)$  (see (\ref{eq: inv}))
 is $\frac{2\pi^2n^2}{\ell}$.  To see this note that in this case
\[
{\en}_\delta (\psi, A_0+a)=\frac{1}{2}\int_N a \wedge da.
\]
Since $F_{A_0+a}^{0,1}=0$ we can pick $a$ of the form $ {\bf i}(t\eta +\pi^*a')$ where $a'\in \Omega^1(\Sigma)$ and $t\in {\bR}$. We get
\[
{\en}_\delta(\psi,A_0+a)=-\frac{t^2}{2}\int_N\eta\wedge d\eta -\frac{t}{2}\int_N
\eta\wedge d a'
\]
Integrating along fibers and using $\int_{N/\Sigma}\eta = 2\pi$ we get
\[
{\en}_\delta(\psi,A_0+a)=\pi t^2\int_\Sigma 2\ell dv_\Sigma  -\pi t \int_\Sigma da' =2\pi^2 t^2 \ell .
\]
On the other hand,  
\[
n=\hat{d}(\psi, A_0+a)=\hat{c}_1(A_0+a)=\frac{\bf i}{2\pi}\int_\Sigma {\bf i}(td\eta +da')=-\frac{-t}{2\pi}\int_\Sigma (-2\ell )dv_\Sigma =\ell t.
\]
Hence
\[
{\en}_\delta(\psi, A_0+a)=\frac{2\pi^2 n^2}{\ell}.
\]
Lemma \ref{lemma: gap} is proved. $\Box$

\bigskip

 To prove the  promised  $L^4_{loc}$ estimates we need the following energy
 identity. 
 
 \begin{lemma}{\rm  If $(\psi(t), a(t))$ is a tunneling of (\ref{eq: aflow})
 then for every $T\in {\bR}$ and any $L>0$ we have
 \[
 \int_{T-L}^{T+L}dt \int_N |{\dps}|^2+ |{\da}|^2 dv_\delta
 \]
 \be
 = \int_{T-L}^{T+L}dt\int_N\left(|\nabla^{A(t)}\psi|^2_\de
 +\frac{1}{8}|\psi|^4+\frac{1}{4}(s_\delta-\lambda^2_\delta)|\psi|^2-\lambda_\delta
 \lan {\dps}, {\psi}\ran +|F_A|_\de^2\right) dv_\delta.
 \label{eq: en4}
 \ee
 Above $A(t)$ denotes the connection $A_0+a(t)$ on the slice $\{t\}\times N$
 and $s_\delta$  denotes the scalar curvature of the metric $h_\delta$.}
\label{lemma: energy}
\end{lemma}

To keep the flow of arguments uninterrupted  we present the proof of this lemma
in an appendix.

We can now produce the promised $L^4_{loc}$ estimates.  Set $C(T, L)=[T-L, T+L]\times
N$. In \cite{Ni} we showed that $\sup_{x\in N} |s_\de(x)|$ is bounded  from above
as $\delta \ra \infty$.  Denote by $|s|$ such an upper bound. The energy identity implies that
there exists a constant $C>0$ independent of  $\delta \gg 1$ such that
\[
\int_{C(T, L)}|\psi|^4dt\,dv_\delta \leq  C\left(E(T, L)+\int_{C(T,L)}|\psi|^2dt\,
dv_\delta\right)
\]
\[
\leq C\left\{E(T,L) + \left(\frac{L}{\delta}\int_{C(T,
L)}|\psi|^4dtdv_\delta\right)^{1/2}\right\}.
\]
The last inequality implies 
\[
\int_{C(T, L)}|\psi|^4dt\,dv_\delta  \leq C(E(T, L)+L\delta^{-1}).
\]
Using the inequality (\ref{eq: en3}) we deduce that there exists a constant $C$ {\em independent of} $\de$ such
that
\be
\int_{C(T, L)}|\psi|^4dt\, dv_\delta \leq C(1+L\delta^{-1}).
\label{eq: en5}
\ee

\subsection{$L^\infty$ estimates}
In the compact case,  effective  $L^\infty$ estimates can be obtained using a
simple maximum principle trick  as in \cite{KM}.  In our noncompact situation
the M\"{o}ser's iteration technique produces results which are dramatically sharper in the adiabatic limit. The starting point is the following Kato type inequality, very
similar to the one used in \cite{KM}. 

\begin{lemma}{\rm  If $(\psi(t), A(t))$ is a solutions of $(\ref{eq:
aflow})_\delta$ then
\be
\Delta_\delta |\psi|^2 +2\lambda_\delta\frac{\partial |\psi|^2}{\partial t} \leq
(\lambda_\delta +z_\delta)|\psi|^2
\label{eq: kato1}
\ee
where $\Delta_\delta$ denotes the scalar Laplacian of the metric
$dt^2+h_\delta$ and $z_\delta=\sup_{x\in N} |s_\delta(x)|$.}
\end{lemma}

\noindent {\bf Proof}\hspace{.3cm}  For simplicity, we prove the above inequality only for
$\delta =1$.  We first rewrite the first equation in (\ref{eq: aflow}) using the operator
$\dol_A$. We get 
\[
\dol_A\psi =  \frac{\lambda}{2}\hat{\bf c}(dt) \psi .
\]
We now apply $\dol_A^*$ to both sides of the above equality  and we obtain
\[
\dol_A^*\dol_A \psi =\frac{\lambda}{2}\dol_A^* \hat{\bf c}(dt).
\]
The anticommutation equality   Prop. 3.45 in \cite{BGV},
\[
{\dol}_A^*\hat{\bf c}(dt)+\hat{\bf c}(dt) {\dol}_A^* =-2\nabla_t
\]
yields
\[
{\dol}_A^*\dol_A \psi= \frac{\lambda^2}{4}\psi -\lambda{\dps}.
\]
Using the Weitzenb\"{o}ck formula in the  left-hand side of the above
inequality we deduce
\[
(\nabla_A^*\nabla_A +\hat{\bf c}(F_A^+) +\frac{s}{4}\psi =\frac{\lambda^2}{4}\psi
-\lambda{\dps}.
\]
The inequality follows  as in \cite{KM},  by taking the (pointwise) inner product with
$\psi$,  then using the equality $\hat{\bf c}(F_A^+)=\tau(\psi)$  to eliminate  a positive  multiple of $|\psi|^4$ and  
concluding with the  Kato inequality (see \cite{Ber})
\[
\Delta|\psi|^2 \leq 2\lan \nabla_A^*\nabla_A\psi, \psi \ran. \;\; \Box
\]

\bigskip

Denote by $S_\delta$  the best Sobolev constant defined  as the largest constant such that
\be
S_\delta \|u\|^2_{4, \de}\leq (\|du\|^2_{2,\de} +\|u\|^2_{2,\de}),\;\;\;\forall
u\in C_0^\infty(C(T,2)).
\label{eq: sobolev}
\ee
Above, $\|\cdot\|_{p,\delta}$ denotes the $L^p$ norm on $C(T,2)$ defined in terms of the
metric $dt^2+h_\delta$. 

We can now state the main result of this subsection.

\begin{lemma}{\rm There exists a constant $C>0$ independent of $\delta$ such
that if $(\psi, A)$ is a tunneling of  $(\ref{eq: aflow})_\delta$ then}
\[
\sup_{(t,x)\in C(T,1)}|\psi(t,x)|^2\leq \frac{C}{S_\delta}\|\psi\|^2_{4,\delta}.
\]
\label{lemma: moser}
\end{lemma}

\noindent{\bf Proof}\hspace{.3cm} Set $u(t,x)=|\psi(t,x)|^2$.   Then $u\in
L^{1,2}_{loc}({\bR}\times N)$ and   satisfies the differential inequality
\be
\Delta_\delta u +2\lambda_\delta \dot{u}\leq  A u \;\;{\rm in}\;\;C(T,2)
\label{eq: mos1}
\ee
where $A$ is a positive constant independent of $\delta$. This is precisely the
type of inequality  for which the M\"{o}ser iteration technique-as described 
e.g. in \cite{GT}- is applicable. This leads to an inequality of the type
\[
\sup_{(t,x)\in C(T,1)}u(t,x)\leq B\|u\|_{2,\delta}.
\]
 To understand the manner in which the constant $B$ depends upon $\delta$ we
 will go carefully through the proof.  To make  the presentation more accessible
  we will omit the subscript $\delta$  in the notations of the various norms
  and in the notation of the volume form. We will also assume $\lambda=0$ in (\ref{eq: mos1})  since this
  term can be  later absorbed anyway via a simple interpolation trick. 
  
  For each $h\in(0,1)$ consider a nonnegative cutoff function 
  \[
  f=f_h(t)\in C_0^\infty(T-2,T+2)
  \]
 identically $1$ on $(T-2+h, T+2-h)$ and such that $|\sup f'(t)|\leq 4h^{-1}$. We set
  for brevity $M=C(T,2)=[T-2,T+2]\times N$ and  $M_h= C(T,2-h)$. 
  
 Multiply the equality
(\ref{eq: mos1}) by $v=f^2u^{2p-1}$ where $p$ will be specified latter. Using
\[
dv=2fu^{2p-1} df + (2p-1)f^2 u^{2p-2}du.
\]
we obtain after an integration by parts
\[
\int_M (2p-1)f^2 u^{2p-2}|du|^2 
\leq A\int_M f^2u^{2p}+\int_M2f u^{2p-1}|du| |df|
\]
\[
=A\int_M (2p-1)f^2 u^{2p-2}|du|^2 +\int_M 2fu^{p-1}|du| \cdot u^p|df|
\]
\[
\leq A\int_M (2p-1)f^2 u^{2p-2}|du|^2 + {\ve}\int_M f^2u^{2p-2}|du|^2 +
{\ve}^{-1}\int_Mu^2 |df|^2.
\]
If we take ${\ve}=1/2$ and   absorb the corresponding term in the
left-hand-side we get
\be
(2p-2)\int_Mf^2 u^{2p-2}|du|^2 \leq  A_1 \int_M u^{2p}
\label{eq: mos2}
\ee
where  $A_1$ is a constant {\em independent} of $\delta$. Now set $w=fu^p$ so
that ${\rm supp}\,w\subset {\rm supp}\, f$ and  on $M$
\[
|dw|^2\leq 2(u^{2p}|df|^2 +p^2 f^2 u^{2p-2}|df|^2).
\]
Hence, using (\ref{eq: mos2}) we deduce
\[
\int_M|dw|^2 +|w|^2 \leq  A_2ph^{-2} \int_M u^{2p}
\]
where $A_2$ is independent of $\delta$.  The Sobolev inequality (\ref{eq:
sobolev}) now implies
\[
S_\delta\|w\|^2_{4,M} \leq A_2ph^{-2}\int_Mu^{2p}.
\]
Hence
\[
\left(\int_{M_h}u^{4p}\right)^{1/2}\leq \frac{A_2p}{S_\delta h^2}\int_{M}u^{2p}
\]
so that
\[
\left(\int_{M_h}u^{4p}\right)^{1/4p}\leq \left(\frac{A_2p}{S_\delta
h^2}\right)^{1/2p}\left(\int_{M}u^{2p}\right)^{1/2p}.
\]
More generally,  for each $0\leq k < h<1$ we have a similar inequality
\be
\left(\int_{M_h}u^{4p}\right)^{1/4p}\leq \left(\frac{zp}{S_\delta
(h-k)^2}\right)^{1/2p}\left(\int_{M_k}u^{2p}\right)^{1/2p}
\label{eq: mos3}
\ee
where $z$ is a positive constant independent of $\delta$, $p$, $h$, $k$. The above inequality is the basis of the M\"{o}ser iterative method.  Set
$p_0=2$, $p_{n+1}=2p_n=2^{n+1}$, $h_0=0$, $h_{n+1}=h_n+2^{-n-1}$. We deduce
\[
\|u\|_{p_{n+1}, M_{h_{n+1}}}\leq
\left(\frac{z2^n}{S_\delta(h_{n+1}-h_{n})^2}\right)^{1/2^n}\|u\|_{p_n,
M_{h_n}}=\left(\frac{z2^{3n+2}}{S_\delta}\right)^{1/2^{n}+1}\|u\|_{p_n,
M_{h_n}}
\]
Iterating and then letting $n\ra \infty$ we deduce
\[
\sup_{M_1}u(x) \leq
\frac{z}{S_\delta}\left(\prod_{n=0}^{\infty}2^{(3n+2)/2^{n+1}}\right)\cdot
\|u\|_{2, M}.
\]
Lemma \ref{lemma: moser} is proved. $\Box$

\subsection{Proof of the Theorem 2.1}
We will argue by contradiction. Assume that for all $\delta \gg 1$ there exists
a tunneling $(\psi(t), A(t))=(\psi_\delta(t), A_\delta(t))$. Since the Ricci curvature of $M$ is bounded as $\delta \ra \infty$ and ${\rm
 diam}\, (M)=4$ we  deduce from  the estimates in \cite{Ber} (formula (2) p. 389
  coupled with Theorem 2 p. 386) that
\[
S_\delta =O({\rm vol}\,(N,h_\delta)^{-1/2})=O(\delta^{1/2})\;\;{\rm as}\; \; \delta \ra \infty.
\]
Lemma \ref{lemma: moser} coupled with (\ref{eq: en5}) yields
\[
\|\psi\|_{\infty, C(T,2)} =O(\delta^{-1/2})\;\;\forall T.
\]
Lemma \ref{lemma: asymptotics} now implies
\be
\|\psi_{\pm}\|_2 =O(\delta^{-1/2})
\label{eq: lim}
\ee
where $\|\cdot\|$ denotes the $L^2$ norm on $N$ with respect to the metric
$h_1$.   At least one of the limits $\psi_\pm$  is  nontrivial since the
reducible part of the critical set is path connected.  On the other hand, by
(\ref{eq: n1}) or (\ref{eq: n2})  we deduce that if $\psi_\pm$ is nontrivial,
then its $L^2$ norm  belongs to a finite set of {\em  universal constants}.
This contradicts (\ref{eq: lim}) and concludes the proof of Theorem \ref{th: vanish}. $\Box$

\bigskip

\begin{remark}{\rm (a) The  result we have proved   suggests that ${\en}_\delta$  is a perfect Morse
function.  To determine the $SW$-Floer (co)homology we must determine the spectral flows of the Hessian
along paths connecting different components of the critical sets which
determine the grading of this homology. We believe it is isomorphic (as a {\em graded} group) with the homology of the
irreducible  part of the adiabatic moduli space. We will deal with this issue elsewhere.
We refer also to \cite{MOY} for an application along these lines.

(b) There  exists a different approach to the study of tunnelings
similar in spirit with  the study of Yang Mills tunnelings in \cite{Guo}. More
precisely,  using (\ref{eq: azt}) and the exponential decay (which is valid  if
 the fibers are sufficiently  short as a consequence of Theorem  \ref{th: bott})
  one can   identify the tunnelings between two
{\em irreducible} components of the critical  set with some  abelian vortices over 
the complex  ruled surface $X$ compactifying ${\bR}\times N$.   One then can
use algebraic geometric techniques to study these vortices. We refer to
\cite{MOY}  where this   approach is dicussed in great detail.}
\end{remark}

\appendix
\section{Proof of Lemma 2.4}
We prove the identity for $\delta=1$. Note that  ${\bf
D}_A^2={\dir}_A^2-\lambda{\dir}_A+\lambda^2/4$. Hence
\[
\int_N|{\dps}|^2 = \int_N|{\bf
D}_A\psi|^2=\int_N\lan {\dir}_A^2\psi, \psi\ran -\lambda\int_N{\re}\lan {\dir}_A\psi, \psi\ran
+\frac{\lambda^2}{4}\int_M|\psi|^2.
\]
Using  the Weitzenb\"{o}ck formula for ${\dir}_A^2$   we deduce
\[
\int_N|{\dps}|^2 =\int_N \left(|\nabla^A\psi|^2+\frac{s(x)}{4}|\psi|^2 +{\re}\lan{\bf
c}(F_A)\psi, \psi\ran\right) -\lambda\int_N{\re}\lan {\dir}_A\psi, \psi\ran
+\frac{\lambda^2}{4}\int_M|\psi|^2.
\]
On the other hand,
\[
2\int_N|\dot{A}|^2=\int_N |{\bf c}(\dot{A})|^2=\int_N|\tau(\psi)|^2 +|{\bf
c}(F_A)|^2 -2\int_N{\re}\lan {\bf c}(F_A), \tau(\psi)\ran
\]
\[
=2\int_N|F_A|^2+\frac{1}{4}\int_N|\psi|^4 -2\int_N{\re}\lan {\bf c}(F_A),
\tau(\psi)\ran.
\]
A simple computation yields
\[
{\re}\lan {\bf c}(F_A), \tau(\psi) \ran ={\re}\lan {\bf c}(F_A)\psi, \psi\ran.
\]
We conclude that
\[
\int_N|{\dps}|^2 +|\dot{A}|^2
\]
\[
=\int_N\left(|\nabla^A\psi|^2+\frac{s(x)}{4}|\psi|^2
+|F_A|^2+\frac{1}{8}|\psi|^4\right) -\lambda\int_N{\re}\lan {\dir}_A\psi, \psi\ran
+\frac{\lambda^2}{4}\int_M|\psi|^2
\]
\[
=\int_N\left(|\nabla^A\psi|^2+\frac{s(x)}{4}|\psi|^2
+|F_A|^2+\frac{1}{8}|\psi|^4\right)-\lambda\int_N{\re}\lan
{\dps}+\frac{\lambda}{2}\psi, \psi \ran +\frac{\lambda^2}{4}\int_N|\psi|^2
\]
\[
=\int_N\left(|\nabla^A\psi|^2+\frac{s(x)}{4}|\psi|^2
+|F_A|^2+\frac{1}{8}|\psi|^4\right)-\frac{\lambda^2}{4}\int_N|\psi|^2-\lambda  
\int_N{\re}\lan{\dps}, \psi \ran.
\] 
This concludes the proof of the energy identity. $\Box$

\end{document}